\documentclass[12pt]{article}

\title{Comparison of arm exponents in planar FK-percolation}

\author{Lo\"ic Gassmann\footnote{D\'epartement de Math\'ematiques,
Universit\'e de Fribourg,
Chemin du Mus\'ee 23,
CH-1700 Fribourg, Switzerland,
Members of NCCR SwissMap,
\href{mailto:loic.gassmann@unifr.ch}{loic.gassmann@unifr.ch},
\href{mailto:ioan.manolescu@unifr.ch}{ioan.manolescu@unifr.ch}} 
\and Ioan Manolescu\footnotemark[1]
}

\usepackage[T1]{fontenc}
\usepackage[english]{babel}
\usepackage{fullpage}
\usepackage{xcolor}  
\usepackage{hyperref}
\usepackage{calrsfs} 
\usepackage{nicefrac}
\usepackage{amsmath}
\usepackage{amsthm}
\usepackage{amssymb}
\usepackage{mathtools}  

\numberwithin{equation}{section}  

\newcommand\Indi[1]{\textbf{1}_{\left\{#1\right\}}}
\newcommand\Ann{{\rm Ann}}

\newcommand\calC{\mathcal{C}}
\newcommand\calF{\mathcal{F}}
\newcommand\calE{\mathcal{E}}

\newcommand\bbN{\mathbb{N}}
\newcommand\bbZ{\mathbb{Z}}
\newcommand\bbP{\mathbb{P}}

\newcommand{\ep}{\varepsilon}
\newcommand\La{\Lambda}

\newtheorem{theorem}{Theorem}[section]
\newtheorem{lemma}[theorem]{Lemma}
\newtheorem{definition}[theorem]{Definition}
\newtheorem{remark}[theorem]{Remark}

\begin{document}

\maketitle

\begin{abstract}
By the FKG inequality for FK-percolation, the probability of the alternating two-arm event is smaller than the product of the probabilities of having a primal arm and a dual arm, respectively. In this paper, we improve this inequality by a polynomial factor for critical planar FK-percolation in the continuous phase transition regime ($1 \leq q \leq 4$). In particular, we prove that if the alternating two-arm exponent~$\alpha_{01}$ and the one-arm exponents~$\alpha_0$ and~$\alpha_1$ exist, then they satisfy the strict inequality~$\alpha_{01} > \alpha_0 + \alpha_1$.

The question was formulated by Garban and Steif in the context of exceptional times and was brought to our attention by Radhakrishnan and Tassion, who obtained the same result for planar Bernoulli percolation through different methods.   
\end{abstract}

\section{Introduction}
A planar percolation model is defined by a probability distribution over the subgraphs of a given planar graph, typically the square lattice~$\bbZ^2$. The most prominent example of these models is Bernoulli percolation on~$\bbZ^2$, where a random subgraph is generated by independently retaining each edge with probability~$p$ or erasing it with probability~$1-p$. On finite graphs, the more general FK-percolation model (also called random-cluster model) is obtained from Bernoulli percolation by adding an exponential bias expressed in terms of the number of connected components, depending on a parameter~$q >0$. The interest of this model comes from its relation to various other statistical physics models, including the Ising and Potts models (see, e.g., \cite{Gri2006}). For~$q=1$ the bias disappears, and we recover the Bernoulli percolation model. One crucial property of Bernoulli percolation, called {\em positive association} is conserved for~$q \geq 1$. This property is expressed by the Fortuin-Kasteleyn-Ginibre (FKG) inequality. 

For~$q \geq1$ fixed, FK-percolation exhibits a phase transition for the existence of an infinite cluster at a certain parameter~$p_c(q)$: 
for~$p < p_c(q)$ there exists almost surely no infinite cluster and for~$p > p_c(q)$ there exists almost surely a unique infinite cluster. In this paper we focus on the critical parameter~$p = p_c(q)$ and restrict ourselves to the regime in which the phase transition is continuous (which is to say~$q \in [1,4]$ as proved in \cite{DcSidTas2017,DcGagHarManTas2016}). 
With this choice of parameters the model possesses a unique infinite Gibbs measure which we denote by~$\phi_{p_c,q}$. This measure is self-dual in  that the primal and the dual configurations have the same law. 
Furthermore, it was proved to exhibit the Russo-Seymour-Welsh (RSW) property, with uniformity in the boundary conditions imposed at a macroscopic distance. As a consequence, certain so-called {\em arm events} involving primal and dual connections occur with polynomially decreasing probabilities. Our goal is to study an inequality between the rates of decay of the probabilities of certain arm-events. We direct the reader to Section~\ref{section:background_FK} for more details on the notions above. 

An arm of type~$1$ (resp. type~$0$) is a path formed of primal open (resp. dual open) edges. We denote by~$\La_n$ the subgraph of~$\bbZ^2$ induced by the vertex set~$\{-n,\dotsb,n\}^2$ and by~$\partial \La_n$ its vertex boundary. For~$0<r<R$, let~$\Ann(r,R)$ denote the annulus~$\La_R \setminus \La_r$. For~$1 \leq r < R$, we denote by~$A_1(r,R)$ (resp.~$A_0(r,R)$) the event that the inner and outer boundaries of~$\Ann(r,R)$ are connected by an arm of type~$1$ (resp. type~$0$). Additionally, denote by~$A_{01}(r,R)$ the event that the inner and outer boundaries of the annulus~$\Ann(r,R)$ are connected by both an arm of type~$0$ and an arm of type~$1$. These events are referred to as arm events. 

For~$1 \leq r < R$, the events~$A_0(r,R)$ and~$A_1(r,R)$ are respectively decreasing and increasing so, by the FKG inequality, we get
\begin{align}\label{ineq:FKG_arms}
	\phi_{p_c,q}[A_{01}(r,R)] \leq \phi_{p_c,q}[A_1(r,R)]\phi_{p_c,q}[A_0(r,R)].
\end{align}
An alternative way of formulating \eqref{ineq:arm_exponent} is that the presence of a dual arm in an annulus hinders that of a primal arm. Our goal is to quantify the effect of the dual arm, namely to prove that it induces an additional cost for the existence of the primal arm which is polynomial in the inner-to-outer radii ratio of the annulus. More precisely, we prove the following.

\begin{theorem}\label{thm:comparison_exponents}
    Fix~$1 \leq q \leq 4$. There exists~$c > 0$ such that for every~$r \leq R$ with~$R/r$ large enough,
    \begin{align}\label{ineq:improved_FKG}
	    \phi_{p_c,q}[A_{01}(r,R)] \leq (r/R)^{c}\, \phi_{p_c,q}[A_1(r,R)]\phi_{p_c,q}[A_0(r,R)].
    \end{align}
\end{theorem}

This inequality can be expressed in terms of arms exponents. Indeed, it is expected that the probabilities of the arm events~$A_{.}(r,R)$ behave as~$r/R$ to some positive exponent, called the one-arm primal, one-arm dual or alternating two-arm exponent, and denoted~$\alpha_1$,~$\alpha_0$ and~$\alpha_{01}$, respectively. Assuming the existence of these exponents, inequality \eqref{ineq:improved_FKG} may be rewritten as
\begin{align}\label{ineq:arm_exponent}
\alpha_{01} > \alpha_0 + \alpha_1.
\end{align}
Note that by the self-duality of~$\phi_{p_c,q}$, we have~$\alpha_0 = \alpha_1$ (if the exponents exist), but we keep the inequality in the form given above for aesthetic reasons.

We mention that it is predicted that the cluster boundaries of critical two-dimensional FK-percolation for~$q \in (0,4]$ converge toward CLE$_{\kappa}$ where 
\begin{align} \label{eq:definition_kappa}
\kappa = 4\pi/\arccos(-\sqrt{q}/2).
\end{align}
Such a convergence would imply the existence of arm exponents and would allow to compute the following expressions for~$q \in (0,4)$ (see \cite{LawSchWer2002,SSW2011} for the one-arm exponent and \cite{SmiWer2001,Wu2018} for the alternating two-arm exponent) 
\begin{align*}
\alpha_0& = \alpha_1 = \frac{(8-\kappa)(3\kappa - 8)}{32 \kappa}  , & 
\alpha_{01}& = \frac{16 - (4-\kappa)^2}{8\kappa},
\end{align*}
where~$\kappa \in (4,8)$ is related to~$q$ via \eqref{eq:definition_kappa}. These values do satisfy~\eqref{ineq:arm_exponent}. Which may appear surprising for~$q \in (0,1)$, as the FKG inequality fails in this regime. So far, the convergence to CLE$_{\kappa}$ was only proven for~$q=2$ \cite{CheSmi2012,CheDcHonKemSmi2014} and for a special case of Bernoulli percolation ($q=1$), namely site percolation on the triangular lattice  \cite{Smi2001,CamNew06}. 

Our proof makes no reference to the supposed scaling limit of the model; it only requires the FKG inequality,  and RSW-type constructions. As such, it is illustrative of a more general phenomenon and provides a general technique to gain polynomial factors in inequalities derived from the FKG inequality. 

As stated previously, Theorem~\ref{thm:comparison_exponents} amounts to comparing the probability of existence of a primal arm in the measure conditioned on the presence of a dual arm and in the unconditioned measure. 
The FKG inequality states that the former is stochastically dominated by the latter, which implies~\eqref{ineq:FKG_arms}. 
It is generally expected that a polynomial improvement may be obtained in such cases, due to differences between the two measures that may appear at every scale with positive probability. Indeed, we generally expect that the existence of a primal arm given a set of existing arms is polynomially more costly than in an unconditioned measure. 

To illustrate this phenomenon, consider critical Bernoulli percolation denoted by $\bbP_{p_c}$ and use the more general notation~$A_\sigma(r,R)$ for arm events with arms of types~$\sigma \in \bigcup_{k \geq 1} \{0,1\}^k$. It is proven in \cite[Appendix A]{GarSte14} that, as long as~$\sigma$ contains at least one~$0$ and one~$1$, there exists a constant~$c > 0$ such that, for every~$r \leq R$ with~$R/r$ large enough,
\begin{align*}
	\bbP_{p_c}[A_{\sigma \cup \{1\}}(r,R)] \leq (r/R)^c \, \bbP_{p_c}[A_\sigma(r,R)] \, \bbP_{p_c}[A_1(r,R)]. 
\end{align*}
This inequality is proven by \emph{exploring} the arms~$\sigma$ and showing that, given the explored regions, the additional primal arm has less space than when no exploration is performed, which results in an additional polynomial cost $(r/R)^c$ for this arm. The difficulty in obtaining \eqref{ineq:improved_FKG} lies in that the dual arm is not explorable without revealing potentially positive information about the existence of a primal arm.

Rather than constructing an appropriate exploration, our proof relies on a carefully designed increasing coupling between ~$\omega \sim \phi_{p_c,q}[\cdot \,|\, A_0(r,R)]$ and~$\omega' \sim \phi_{p_c,q}[\cdot]$. We show that in this coupling,~$\omega'$ is polynomially more likely to contain a primal arm than~$\omega$. To do this, we construct a special series of scales (called {\em good scales}) which each contribute a multiplicative constant to the difference between the arm probabilities in~$\omega$ and~$\omega'$. 

\begin{remark}
One can adapt the proof of Theorem~\ref{thm:comparison_exponents} to show that inequality~\eqref{ineq:improved_FKG} remains valid if we replace $A_{01}(r,R)$ by $A_{\sigma 1}(r,R)$ and $A_{0}(r,R)$ by $A_{\sigma}(r,R)$, where $\sigma$ is any finite sequence of zeros. We do not provide the details of this adaptation, except to note that it requires to change the definition of good scales by including additional boxes containing flower domains. Pivotal edges are then employed to ensure that some boxes cannot be crossed by more than one arm of the annulus.
\end{remark}

The result of Theorem~\ref{thm:comparison_exponents} was conjectured for two-dimensional critical Bernoulli percolation in \cite[Open problem XIII.6]{GarSte14} as a shorter path towards proving the existence of exceptional times (see also \cite{GPS10}). In this setting, the result was proved by Radhakrishnan and Tassion \cite{RadTas24} through different methods. Indeed, rather than constructing an increasing coupling, they interpolate between~$\omega \sim \phi_{p_c,1}[\cdot \,|\, A_0(r,R)]$ and~$\omega' \sim \phi_{p_c,1}$ in a dynamical fashion. While their technique only applies to Bernoulli percolation, they provide more general inequalities, in particular a quantitative improvement of Reimer's inequality reflected in a difference between the monochromatic and alternating two-arm events.\\

\noindent\textbf{Organisation of the paper:} In section~\ref{section:background_FK} we provide some minimal background on FK-percolation, with references for further reading. Section~\ref{section:scale_construction} contains the definition of good scale and proves the existence of a density of such scales in the annulus~$\Ann(r,R)$. In Section~\ref{section:increasing_coupling}, we present the increasing coupling between the conditioned and the unconditioned measures and we prove that a density of good scales implies an additional polynomial cost for the existence of the primal arm in~$\omega$ than~$\omega'$. Finally, in Section~\ref{section:proof_main}, we combine the results of the two previous sections to prove Theorem~\ref{thm:comparison_exponents}.

\paragraph{Acknowledgements:} We thank Ritvik Ramanan Radhakrishnan and Vincent Tassion for bringing this problem to our attention and encouraging us to publish this alternative proof of their result \cite{RadTas24}. We also thank Christophe Garban for useful discussions on the topic. 

Part of this project was developed during the first author's internship at the University of Fribourg in 2024. During this internship the first author was supported by the University of Fribourg and the \'Ecole Normale Sup\'erieure Paris-Saclay. The authors are part of the NCCR SwissMAP.

\section{Background on FK-percolation}\label{section:background_FK}

This section contains a brief introduction to FK-percolation. For a more comprehensive exposition, we refer to the monograph \cite{Gri2006} or the more recent lecture notes \cite{DC2020}. Consider a finite subgraph~$G$ of~$\bbZ^2$ with vertex set $V$ and edge set $E$. The FK-percolation on the graph $G$ is a random measure on the set~$\{0,1\}^E$ of percolation configurations. In a fixed configuration~$\omega \in \{0,1\}^E$, an edge~$e$ of~$E$ is said to be open if~$\omega_e = 1$ and closed otherwise. The configuration~$\omega$ is identified with the set of open edges as well as with the graph constituted of the vertices of~$G$ and the open edges of~$\omega$. An open (or primal) path in a percolation configuration is a sequence of vertices such that each consecutive pair of vertices in the sequence is connected by an open edge. We say that two sets of vertices~$V_1$ and~$V_2$ are connected by an open path if there exists an open path that starts at a vertex of~$V_1$ and ends at a vertex of~$V_2$. A cluster of a configuration~$\omega$ is a connected component of the graph~$\omega$.

The boundary~$\partial G$ of~$G$ is the set of vertices of~$G$ with neighbours in~$\bbZ^2 \setminus G$. 
A boundary condition~$\xi$ on~$G$ is a partition of the vertices of~$\partial G$. Points on the boundary that belong to the same element of~$\xi$ are said to be wired together. We denote by~$\omega^\xi$ the graph obtained from~$\omega$ by identifying all the vertices that are wired together. When all the boundary points are wired together, we obtain the wired boundary conditions (denoted by~$1$), and when no points are wired together, we obtain the free boundary conditions (denoted by~$0$). 

\begin{definition}
The FK-percolation on~$G$ with edge weight~$p \in (0,1)$, cluster weight~$q > 0$ and boundary conditions~$\xi$ is a probability measure on the percolation configurations~$\omega~\in~\{0,1\}^E$, defined by
\begin{align*}
\phi_{G,p,q}^\xi[\omega]=\frac{1}{Z^\xi(G,p,q)} \Big(\frac{p}{1-p}\Big)^{|\omega|}q^{k(\omega^\xi)},
\end{align*}
where~$|\omega|$ denotes the number of open edges in~$\omega$,~$k(\omega^\xi)$ denotes the number of connected components of~$\omega^\xi$ and~$Z^\xi(G,p,q)$ is a normalising constant ensuring that~$\phi_{G,p,q}^\xi$ is indeed a probability measure.
\end{definition}

For~$q \geq 1$, the measures~$\phi_{G,p,q}^0$ and~$\phi_{G,p,q}^1$ converge weakly as~$G$ tends to the entire square lattice. We denote by~$\phi_{p,q}^0$ and~$\phi_{p,q}^1$ the corresponding infinite volume measures.
In this paper we are interested in the critical point
\begin{align*}
p_c = p_c(q) = \inf\{p \in (0,1) \,|\, \phi_{p,q}^1[0 \text{ is in an infinite cluster}] > 0\}.
\end{align*}
The critical probability~$p_c$ was proven in \cite{BefDC2012} to be equal to the self-dual point of the model:
\begin{align*}
	p_c(q) = \frac{\sqrt{q}}{1+\sqrt{q}}.
\end{align*}

For~$q \in [1,4]$ the phase transition is continuous \cite{DcSidTas2017}, which means that~$\phi_{p_c,q}^0 = \phi_{p_c,q}^1$ and $\phi_{p_c,q}^1[0 \text{ is in an infinite cluster}] = 0$. In this case, we simply write~$\phi_{p_c,q}$ for the critical infinite-volume measure. For~$q >4$ the phase transition is discontinuous \cite{DcGagHarManTas2016}, and~$\phi_{p_c,q}^0$ and~$\phi_{p_c,q}^1$ are distinct and have sub- and super-critical behaviours, respectively. 
Our interest lies in critical phenomena, and we will focus on~$q \in [1,4]$ and~$p = p_c(q)$.

Below is a list of additional properties of the model that we will use in this paper. These properties are stated briefly and some of them are specific to particular parameter choices. We refer the reader to \cite{Gri2006,DC2020} for a broader presentation of these properties.\\

\noindent\textbf{Monotonic properties.} We say that an event~$A$ is increasing if, for all~$\omega \leq \omega'$,~$\omega \in A$ implies that~$\omega' \in A$. Fix~$q \geq 1$,~$p \in (0,1)$ and a subgraph~$G$ of the square lattice. For~$A$ and~$B$ two increasing events, the FKG inequality states that~$A$ and~$B$ are positively correlated, i.e.
\begin{align}\label{eq:FKG}\tag{FKG}
\phi_{G,p,q}^{\xi}[A \cap B] &\geq \phi_{G,p,q}^{\xi}[A] \phi_{G,p,q}^{\xi}[B].
\end{align}
For~$\xi$ and~$\xi'$ two boundary conditions on~$G$ such that~$\xi \leq \xi'$, which means that any pair of vertices wired in~$\xi$ is also wired in~$\xi'$, the FK-percolation with boundary condition~$\xi'$ dominates the one with boundary conditions~$\xi$. This means that for all increasing event~$A$,
\begin{align}\label{eq:Mon}\tag{MON}
\phi_{G,p,q}^{\xi'}[A] &\geq \phi_{G,p,q}^{\xi}[A].
\end{align}

\noindent\textbf{The spatial Markov property.} For~$\omega'$ a percolation configuration on~$\bbZ^2$ and~$F$ a finite subset of the edge set of~$\bbZ^2$, we have
\begin{align}\label{eq:SMP}\tag{SMP}
\phi_{p_c,q}[\cdot_{|F} \,|\, \omega_e=\omega'_e,\forall e\notin F] = \phi^\xi_{H,p_c,q} [\cdot],
\end{align}
where~$H$ is the graph induced by the edge set~$F$ and~$\xi$ is the boundary condition induced on~$H$ by the restriction of~$\omega'$ to the complement of~$F$. This means that two vertices of~$\partial H$ are wired in~$\xi$ if and only if there exists an open path of~$\omega'$ between them in the complement of~$F$.\\

\noindent\textbf{Dual model.} The dual graph associated with~$\bbZ^2$ is the graph~$(\bbZ^2)^*$ constructed by placing a vertex at the centre of each square of~$\bbZ^2$ and, for each edge~$e$ of~$\bbZ^2$, adding an edge~$e^*$ between the two vertices corresponding to the faces bordering~$e$. To each configuration~$\omega$ on~$\bbZ^2$, we associate the dual configuration~$\omega^*$ on~$\bbZ^2$ defined by taking~$\omega^*_{e^*} = 1 - \omega_e$, for all edges~$e$ of~$\bbZ^2$. A dual path is a path constituted of edges~$e^*$ such that~$\omega^*_{e^*} = 1$. 
The self-duality of~$p_c$ and the equality~$\phi_{p_c,q}^0 = \phi_{p_c,q}^1$ imply that if~$\omega \sim  \phi_{p_c,q}$, then~$\omega^*$ and~$\omega$ have the same law (up to a shift by~$(1/2,1/2)$).\\

\noindent\textbf{Russo-Seymour-Welsh (RSW) theory.} 
For~$q \in [1,4]$ and~$p = p_c(q)$, the probabilities of crossing rectangles of a given aspect ratio but arbitrary scale is uniformly bounded away from~$0$ and~$1$. Moreover, these bounds remain uniform even when boundary conditions are imposed at a macroscopic distance from the rectangle. Specifically, for~$\rho, \ep > 0$, there exists~$c > 0$ such that for every graph~$G$ containing the rectangle~$[-\ep n,(\rho+\ep)n] \times [-\ep n, (1+\ep)n]$ and every boundary conditions~$\xi$,
\begin{align}\label{eq:RSW}\tag{RSW}
c \leq \phi_{G,p_c,q}^{\xi}[\calC([0,\rho n] \times [0,n])] \leq 1-c,
\end{align}
where~$\calC([0,\rho n] \times [0,n])$ is the event that the rectangle~$[0,\rho n] \times [0,n]$ contains an open path between its left and right sides (see \cite{DcSidTas2017}).

One classical consequence of the RSW theory, combined with (FKG) and (SMP), is that arm events satisfy a quasi-multiplicativity property and arm-separation. It also implies that the probabilities of the arm events $A_0(r,R)$, $A_1(r,R)$ and $A_{01}(r,R)$ may be bounded above and below by polynomials of~$r/R$, with strictly positive powers. Note that this last property is weaker than the existence of arm exponents, which is conjectured by conformal invariance theory. 

In the present paper we use a form of arm-separation for the alternating two-arm event~$A_{01}(r,R)$, which will be stated precisely later on. 
The theory of arm-separation has been developed for Bernoulli percolation in \cite{Kes87} and reworked in \cite{Nol2008}. 
Due to \eqref{eq:RSW}, it may be adapted to alternating arm events in critical FK-percolation with~$q \in [1,4]$.
To avoid the use of ``locally monotone events'' as in \cite[Lem.~13]{Nol2008}, it is better to work with the notion of well-separated flower domains (defined below) instead of the well-separated scales of \cite{Nol2008}. 
Indeed, flower domains may be explored from one side and, if well-separated, connections to their petals may be constructed by direct applications of \eqref{eq:RSW}, \eqref{eq:SMP} and \eqref{eq:Mon} --- see the proof of Lemma~\ref{lem:density_good_scales} for an example. 

With these minor modifications, the proof of arm-separation of \cite{Nol2008} adapts readily to FK-percolation, and we give no further details. 
For a proof of arm-separation for more general arm events, and which relies on the stronger RSW inequalities available for $q \in [1,4)$, the interested reader may look at \cite[Section 5]{CheDcHon2016} (see also \cite{DcManTas2021}). 
We also refer to \cite{DM2022} for examples of the use of arm-separation in FK-percolation.\\

\noindent\textbf{Increasing coupling via exploration.} Consider a decreasing event~$A$ (i.e. an event whose complementary is increasing). Then \eqref{eq:FKG} states that there exists a coupling~$\bbP$ between~$\omega \sim \phi[\cdot \,|\, A]$ and~$\omega' \sim \phi[\cdot]$ such that~$\omega \leq \omega'$ a.s.

It will be useful to construct an explicit such coupling as follows. Fix an ordering~$e_1,e_2,\dots$ of the edges and sample each~$\omega(e_i),\omega'(e_i)$ sequentially according to their respective laws conditioned on the state of the previously sampled edges. If we assume that~$\omega(e_j) \leq \omega'(e_j)$ for all $j < i$, then \eqref{eq:FKG}, \eqref{eq:SMP} and \eqref{eq:Mon} imply that the law of~$\omega(e_i)$ is dominated by that of~$\omega'(e_i)$. Thus, one may sample the state of $e_i$ such that~$\omega(e_i) \leq \omega'(e_i)$ a.s. By induction, it follows that~$\omega \leq \omega'$. 

Notice that the ``exploration'' order $e_1,e_2,\dots$ may also be adapted during the sampling, with the choice of $e_i$ depending on the states of the sampled edges $e_j$ for $j<i$. See \cite[Sec 2.3]{DM2022} for more details.\\

\begin{figure}
  \centering
  \includegraphics{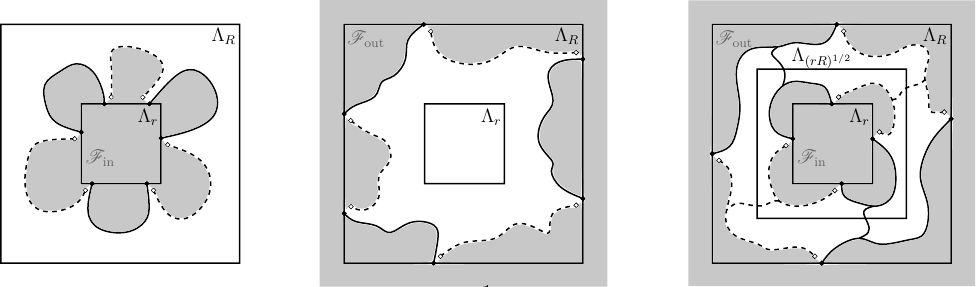}
\caption{An inner flower domain, an outer flower domain and a double four-petal flower domain. The plain and dashed lines represent primal and dual paths respectively.} 
\label{fig:flower_domains}
\end{figure}

\noindent\textbf{Flower domains.} Given~$r \leq R$ and a configuration~$\omega$, we define the inner flower domain $\calF_{\rm in}$ between~$\Lambda_R$ and~$\Lambda_r$ as follows. Consider all interfaces between primal and dual clusters of~$\omega$ and $\omega^*$, respectively, contained in~$\La_R \setminus \La_r$ and starting on~$\partial \La_R$. Write~$\calF_{\rm in}$ for the connected component of the complement of this set of interfaces that contains $\La_r$. In a similar fashion, we define the outer flower domain~$\calF_{\rm out}$ between~$\La_r$ and~$\La_R$ by considering the interfaces in the annulus $\La_R \setminus \La_r$ starting on~$\partial\La_r$, and setting $\calF_{\rm out}$ to be the infinite connected component of the complement of this set of interfaces.  See Figure~\ref{fig:flower_domains} for an illustration and \cite{DM2022} for more detailed definitions.
A crucial feature of the flower domains $\calF_{\rm in}$ and $\calF_{\rm out}$ is that they are measurable in terms of the edges 
of $\La_R \setminus \calF_{\rm in}$ and $(\La_r \cup \calF_{\rm out})^c$, respectively. 

The boundaries of inner and outer flower domains are naturally split in arcs formed alternatively by the primal and dual sides of the interfaces used to define them; we call these arcs the {\em petals} of the flower domain.

If none of the interfaces starting at~$\La_R$ touches~$\partial \La_r$, then~$\calF_{\rm in}$ has only one petal, which is a primal or dual circuit around~$\La_r$. In all other cases the number of petals of $\calF_{\rm in}$ is even and their endpoints all lie on $\partial \La_R$. 
Similar considerations hold for $\calF_{\rm out}$. 
The flower domain~$\calF_{\rm in}$ (resp.~$\calF_{\rm out}$) is said to be well-separated if it contains at least two petals and the endpoints of each petal are separated by a distance greater than~$r/2$ (resp.~$R/2$).

We say that there exists a double four-petal flower domain between~$\La_r$ and~$\La_R$ if the following condition are met. The inner flower domain~$\calF_{\rm in}$ between~$\La_{(rR)^{1/2}}$ and~$\La_r$ and the outer flower domain~$\calF_{\rm out}$ between~$\La_{(rR)^{1/2}}$ and~$\La_R$ have four petals each and are well-separated. We also require that each primal (resp. dual) petal of~$\calF_{\rm in}$ is connected by a primal (resp. dual) path to a primal (resp. dual) petal of~$\calF_{\rm out}$, as illustrated in Figure~\ref{fig:flower_domains}. The following lemma \cite[Lemma 3.4]{DM2022} allows us to construct such a double four-petal flower domain with positive probability. 

\begin{lemma}[\cite{DM2022}]\label{lem:four_petal_flower_domain}	
	For any $q \in [1,4]$ and~$\eta > 0$ there exists $c > 0$ such that the following holds. For any $R$ large enough, and any boundary conditions $\xi$ on $\La_{(1+2\eta) R}$, 
	\begin{align*}
	\phi_{\La_{(1+2\eta) R}, p_c,q}^\xi \big[\text{there exists a double four-petal flower domain between~$\La_R$ and~$\La_{(1+\eta)R}$}\big] > c. 
	\end{align*}
\end{lemma}

\medskip

\begin{center}\bf
	For the rest of the paper, we fix~$q \in [1,4]$,~$p=p_c(q)$ and we write~$\phi$ instead of~$\phi_{p_c,q}$.
\end{center}

\section{Density of good scales}\label{section:scale_construction}

In preparation for the construction of the coupling that we will use to prove Theorem~\ref{thm:comparison_exponents}, we define a notion of {\em good scale}, and show that a linear number of such scales appear between $r$ and $R$ for the measure conditioned on $A_{01}(r,R)$. 
For now we work with a single configuration. 

\begin{figure}[t]
  \centering
  \includegraphics{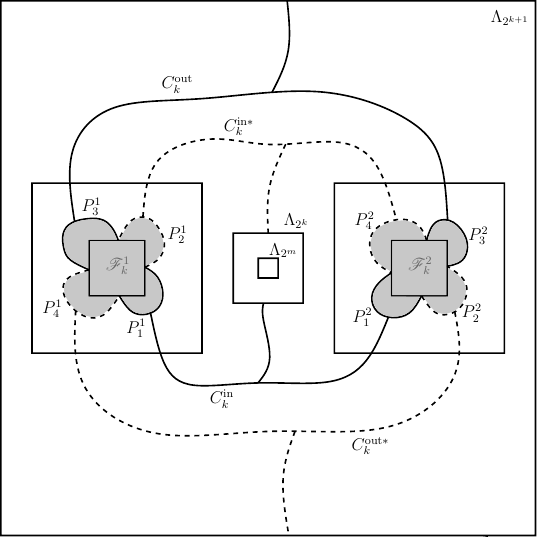}
	\caption{A depiction of a configuration in which the scale~$k$ is good.} \label{fig:good_scale}
\end{figure}

For simplicity, we will work with $r = 2^m$ and $R = 2^n$ with~$m,n \in \bbN$. For~$m \leq k < n$, we will refer to the annulus~$\Ann(2^k,2^{k+1})$ as the scale~$k$. For a fixed $k$, define the following square boxes contained within the scale~$k$
\begin{align*}
S^1_k &= \La_{\frac{2^k}{10}} + 2^k(-\tfrac{3}{2},0)  , &
S^2_k &= \La_{\frac{2^k}{10}} + 2^k(\tfrac{3}{2},0),\\ 
\hat{S}^1_k &= \La_{\frac{2^k}{5}} + 2^k(-\tfrac{3}{2},0)  ,&
\hat{S}^2_k &= \La_{\frac{2^k}{5}} + 2^k(\tfrac{3}{2},0).
\end{align*}

\begin{definition}
We say that the scale~$k$ is good (see Figure~\ref{fig:good_scale}) if
the inner flower domains~$\calF_k^{i}$ between~$\hat{S}^i_k$ and~$S^i_k$ are well-separated and contain four petals~$P_1^i, \dots, P_4^i$ for~$i=1,2$ and the following connections occur in~$\Ann(2^k,2^{k+1}) \setminus (\calF_k^{1} \cup \calF_k^{2})$:

\begin{itemize}
\item $P_1^1$ and~$P_1^2$ are connected by a primal cluster~$C^{\rm in}_k$ that also intersects~$\partial \La_{2^k}$;
\item $P_3^1$ and~$P_3^2$ are connected by a primal cluster~$C^{\rm out}_k$ that also intersects~$\partial \La_{2^{k+1}}$;
\item $P_2^1$ and~$P_4^2$ are connected by a dual cluster~$C^{\rm in *}_k$ that also intersects~$\partial \La_{2^k}$;
\item $P_4^1$ and~$P_2^2$ are connected by a dual cluster~$C^{\rm out *}_k$ that also intersects~$\partial \La_{2^{k+1}}$.
\end{itemize}
We say that the good scales~$k_1<\dots< k_\ell$ are in series if the following connections occur in $\Ann(2^m,2^{n}) \setminus \bigcup_{i,j}\calF_{k_j}^{i}$ (with the union over~$i=1,2$ and~$1\leq j \leq \ell$):
\begin{itemize}
\item for each~$1 \leq j < \ell$,~$C^{\rm out}_{k_j}$ is connected to~$C^{\rm in}_{k_{j+1}}$ by a primal path;
\item for each~$1 \leq j < \ell$,~$C^{\rm out *}_{k_j}$ is connected to~$C^{\rm in *}_{k_{j+1}}$ by a dual path;
\item $C^{\rm in}_{k_1}$ is connected to~$\partial \La_{2^m}$ and~$C^{\rm out}_{k_\ell}$ to~$\partial \La_{2^n}$ by primal paths;
\item $C^{\rm in *}_{k_1}$ is connected to~$\partial \La_{2^m}$ and~$C^{\rm out *}_{k_\ell}$ to~$\partial \La_{2^n}$ by dual paths.
\end{itemize}
\end{definition}

Denote by~$G_k$ the event that the scale~$k$ is good. Additionally, for~$K>0$ we write~$D(K) = D_{m,n}(K)$ the event that there are at least~$K$ good scales in~$\Ann(2^m,2^n)$ which are in series. 
Notice that  $D(K)$, as well as the set of good scales in series $\{k_1,k_2,\dots,k_\ell\}$ are measurable in terms of the edges outside of $\bigcup_j(\calF_{k_j}^{1} \cup \calF_{k_j}^{2})$; 
an explicit exploration will be descried in Section~\ref{section:increasing_coupling}.

\begin{lemma}[Positive density of good scales]\label{lem:density_good_scales}
    There exists~$\delta > 0$ such that, for all~$n\geq m \geq 1$ large enough,
    \begin{align}\label{eq:density_good_scales}
    \phi\Big[D(\delta (n-m)) \,\Big|\, A_{01}(2^m,2^n)\Big] \geq \delta.
    \end{align}
\end{lemma}

The experienced reader will notice that the above is a standard consequence of arm-separation.
Moreover, the probability in \eqref{eq:density_good_scales} may actually be shown to be greater than~$1-e^{-\delta(n-m)}$ for some potentially altered value of~$\delta$. We stated the weaker bound \eqref{eq:density_good_scales} as this suffices for our purposes and may be derived through a simple first-moment argument. 

\begin{proof}[Proof of Lemma~\ref{lem:density_good_scales}]
%It suffices to prove the Lemma for~$n-m$ sufficiently large. \gl{I don't see where we could need that $n-m$ is large in the sequel}. 
Observe that, if~$k_1,\dots,k_\ell$ are good scales and~$A_{01}(2^m,2^n)$ occurs, then they are necessarily in series. 
Thus, our goal is to prove that, under the conditioning~$A_{01}(2^m,2^n)$ the number of good scales is linear with positive probability. 
We will use the term uniform constant to mean a constant $c > 0$ that is independent of $n$, $m$ or $k$. 

We will use a first-moment method, and start by estimating the probability that a fixed scale is good. 
Thus, we begin by proving that for~$m \leq k < n$,
\begin{align}\label{ineq:positive_G_conditional}
	\phi[G_k \,|\, A_{01}(2^m,2^n)] \geq c,
\end{align}
for some uniform constant~$c >0$. 

The proof of this bound consists in constructing the good event~$G_k$ step by step using an exploration and proving that each step of the construction has a positive probability of success. In the following we fix~$m+1 \leq k < n-1$ (for simplicity's sake, we do not treat the cases~$k=m$ and~$k=n-1$).

Start by exploring the outer flower domain~$\calF$ from~$\La_{2^{k-1}}$ to~$\La_{2^k}$ and the inner flower domain~$\calF'$ from~$\La_{2^{k+2}}$ to~$\La_{2^{k+1}}$, as well as the annuli $\Ann(2^m,2^{k-1})$ and $\Ann(2^{k+2},2^n)$. Then, for~$i=1,2$ explore the double four-petal flower domains between~$S^i_k$ and~$\hat{S}^i_k$. Denote by~$\calF_{\rm in}^i$ (resp.~$\calF_{\rm out}^i$) the corresponding inner (resp. outer) flower domain. Write~$E$ for the set of edges revealed at this time (see Figure~\ref{fig:good_scale_construction}). Call~$H$ the event that we obtain two well-separated flower domains at the first step and~$H'$ the event that we obtain two double four-petal flower domains at the second step. With these notations, it follows from the separation of arms for the alternating two-arm event that
\begin{align*}
\phi\big[H \,\big|\, A_{01}(2^m,2^{k}) \cap A_{01}(2^{k+1},2^{n})\big] \geq c_1,
\end{align*} 
where~$c_1 > 0$ is a uniform constant. Then, Lemma~\ref{lem:four_petal_flower_domain} yields a uniform constant~$c_2 > 0$, such that
\begin{align}\label{ineq:positive_H'}
\phi\big[H \cap H' \,\big|\, A_{01}(2^m,2^{k}) \cap A_{01}(2^{k+1},2^{n})\big] \geq c_2.
\end{align}

\begin{figure}[t]
  \centering
  \includegraphics{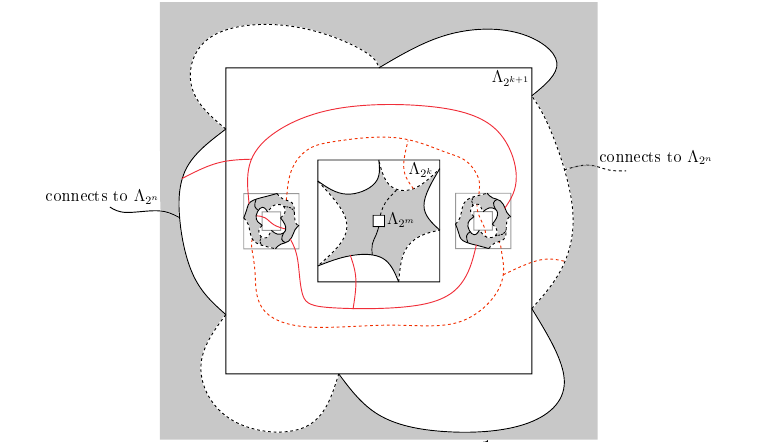}
\caption{In grey the explored region~$E$. {If $H \cap H' \cap A_{01}(2^m,2^{k}) \cap A_{01}(2^{k+1},2^{n})$ occurs in $E$}, then the additional red connections outside of $E$ are sufficient to ensure that~$G_k \cap A_{01}(2^m,2^{n})$ occurs.} \label{fig:good_scale_construction}
\end{figure}

We will now consider the measure $\phi$ conditionally on the configuration $\omega$ on $E$ satisfying the event $H \cap H' \cap A_{01}(2^m,2^{k}) \cap A_{01}(2^{k+1},2^{n})$. Note that all the events of the intersection are indeed measurable in terms of $\omega$ on $E$. 
Then, for the event~$G_k \cap A_{01}(2^m,2^{n})$ to be realised, it suffices that the following connections occur in~$\Ann(2^{k-1},2^{k+2}) \setminus E$ (see the red connections in Figure~\ref{fig:good_scale_construction}):
\begin{itemize}
\item a primal connection between the two primal petals of~$\calF_{\rm in}^1$, contained in~$\calF_{\rm in}^1$,
\item a dual connection between the two dual petals of~$\calF_{\rm in}^2$, contained in~$\calF_{\rm in}^2$,
\item some primal and dual connections between the petals of~$\calF$,~$\calF'$,~$\calF_{\rm out}^1$ and~$\calF_{\rm out}^2$, as specified in Figure~\ref{fig:good_scale_construction}.
\end{itemize}
The three points above only require connections between petals of well-separated flower domains to occur in the unexplored region~$E^c$. So, by a standard application of \eqref{eq:RSW}, \eqref{eq:SMP} and \eqref{eq:Mon}  we get, for a uniform constant~$c_3 > 0$,
\begin{align*}
	\phi\big[G_k \cap A_{01}(2^m,2^n) \,\big|\, \omega \text{ on }E\big] \geq c_3,
\end{align*} 
whenever $\omega$ on $E$ is such that $H \cap H' \cap A_{01}(2^m,2^{k}) \cap A_{01}(2^{k+1},2^{n})$ occurs.
Averaging over the realisations of $\omega$ on $E$ we find
\begin{align}\label{ineq:positive_G}
\phi\big[G_k \cap A_{01}(2^m,2^n) \,\big|\,  H \cap H' \cap A_{01}(2^m,2^{k}) \cap A_{01}(2^{k+1},2^{n})\big] \geq c_3.
\end{align} 
Combining inequalities \eqref{ineq:positive_H'} and \eqref{ineq:positive_G}, yields
\begin{align*}
\phi\big[G_k \cap A_{01}(2^m,2^{n}) \,\big|\, A_{01}(2^m,2^{k}) \cap A_{01}(2^{k+1},2^{n}) \big] \geq c_2c_3. 
\end{align*}
By inclusion of events, we conclude that \eqref{ineq:positive_G_conditional} holds for some uniform constant~$c > 0$. Hence
\begin{align*}
\phi\Big[\sum\limits_{k=m}^{n-1} \Indi{G_k} \,\Big|\, A_{01}(2^m,2^n)\Big]  \geq c (n-m),
\end{align*}
which yields, for all~$\delta > 0$,
\begin{align*}
\phi\Big[\sum\limits_{k=m}^{n-1} \Indi{G_k} \geq \delta (n-m) \,\Big|\, A_{01}(2^m,2^n)\Big] \geq c-\delta.
\end{align*}
Choosing~$\delta = c/2$ concludes the proof.
\end{proof}

\section{An increasing coupling}\label{section:increasing_coupling}

Below,~$\bbP$ will be an increasing coupling of~$\omega \sim \phi[\cdot \,|\, A_0(2^m,2^n)]$ and~$\omega' \sim \phi[\cdot]$, i.e. a coupling satisfying that~$\omega \leq  \omega'$ a.s. Since the event~$A_0(2^m,2^n)$ is decreasing, we know from Section~\ref{section:background_FK} that such coupling exists and may be constructed by exploration. 

\begin{lemma}\label{lem:comparison_under_good_event}
    There exists~$C > 1$ such that, for all~$m \leq n$ and~$K \leq n-m$, there exists a coupling~$\bbP$ as above such that 
    \begin{align*}
    \bbP\big[\omega \in A_1(2^m,2^n) \,\big|\, \omega \in D(K)\big] \leq \big(\tfrac{1}{C}\big)^K \,\bbP\big[\omega' \in A_1(2^m,2^n) \,\big|\, \omega \in D(K)\big].
    \end{align*}
\end{lemma}

\begin{proof}
Fix integers~$m < n$ and $K \leq n-m$. We start by constructing the coupling~$\bbP$. We do so by revealing the state of edges in the two configurations~$\omega$ and~$\omega'$ in a specific order that depends on the configuration~$\omega$. Each time we reveal the state of an edge in $\omega$, we also reveal the state of the same edge in~$\omega'$.

For each~$m\leq k<n$, start by revealing the edges of~$\Ann(2^k,2^{k+1}) \setminus (\hat{S}^{1}_k \cup \hat{S}^{2}_k)$.
Then, explore the outside of the inner flower domains~$\calF^{i}_k$ in~$\omega$ between~$\hat{S}^{i}_k$ and~${S}^{i}_k$ for~$i=1,2$.
At this stage, we may decide if the scale~$k$ is good in~$\omega$; if this is not the case, reveal the entire configuration in~$\Ann(2^k,2^{k+1})$.

After doing this procedure at each scale~$k$, call~$k_1,\dots,k_\ell$ the good scales. Notice that we may now determine if they are in series, and therefore if~$D(K)$ occurs for~$\omega$. If this is not the case, reveal the rest of the configuration. 

Assume henceforth that~$\omega \in D(K)$. Write~$E$ for the set of edges revealed at this time: 
\begin{align}
	E =E(\omega) = \Ann(2^m,2^n) \setminus \bigcup_{j=1}^{\ell} (\calF^{1}_{k_j} \cup  \calF^{2}_{k_j}).
\end{align}
We will now sample the configurations~$\omega$ and~$\omega'$ inside each~$\calF^{i}_{k_j}$. 
Fix some realisations~$\omega_0,\omega_0'$ such that~$\omega_0 \in D(K)$ and write~$\calE$ for the event that~$(\omega,\omega') = (\omega_0,\omega'_0)$ on~$E$. 

For clarity, we will sample a third configuration~$\tilde \omega$ according to~$\phi[ . \,|\, \omega \text{ on } E]$.
Note that,~$\tilde \omega \leq_{\rm st}  \omega'$, since~$\omega \leq \omega'$ on~$E$, and therefore the boundary conditions induced on the domains~$\calF^{i}_{k_j}$ by~$\omega'$ are greater than those induced by~$\omega$.
Additionally~$\omega \leq_{\rm st} \tilde \omega$, as the two configurations have the same boundary conditions, but~$\omega$ is conditioned on the decreasing event~$A_0(2^m,2^n)$. 
Sample these three configurations so that~$\omega \leq \tilde\omega \leq \omega'$ and~$\tilde\omega = \omega$ on~$E$. 

For~$i=1,2$ and~$1 \leq j\leq \ell$ we say that~$\calF_{k_j}^i$ is primally crossed if its two primal petals are connected in~$\calF_{k_j}^i$. Otherwise we say that it is dually crossed.  Write~$X_j^i$ and~$Y_j^i$ for the indicators that~$\calF_{k_j}^i$ is primally crossed in~$\omega$ and $\tilde \omega$, respectively. 

Notice that, conditionally on~$\calE$, the couples~$(Y_j^1 ,Y_j^2)$ are independent for different values of~$j$. Indeed, the realisations of~$\tilde\omega$ inside the flower domains at one scale do not influence the boundary conditions at a different scale (see Figure~\ref{fig:density_good_events}). For a fixed $j$,~$Y_j^1$ and~$Y_j^2$ are not independent of each other, but the RSW theory allows to show that the couple~$(Y_j^1 ,Y_j^2)$ takes each of the four possibilities~$\{0,1\}^2$ with uniformly positive probability. 

The requirement that~$\omega \in A_0(2^m,2^n)$ is equivalent to having at least one of~$\calF_{k_j}^1$ or~$\calF_{k_j}^2$ being dually crossed at each scale. In other words the conditioning imposes that~$X_j^i + X_j^i \leq 1$ for each~$j$. 
Thus, the law of the variables~$(X_j^i)_{i,j}$ is that of~$(Y_j^i)_{i,j}$ conditionally on~$Y_j^i + Y_j^i \leq 1$ for each~$j$.

\begin{figure}[t]
  \centering
  \includegraphics[scale=0.9]{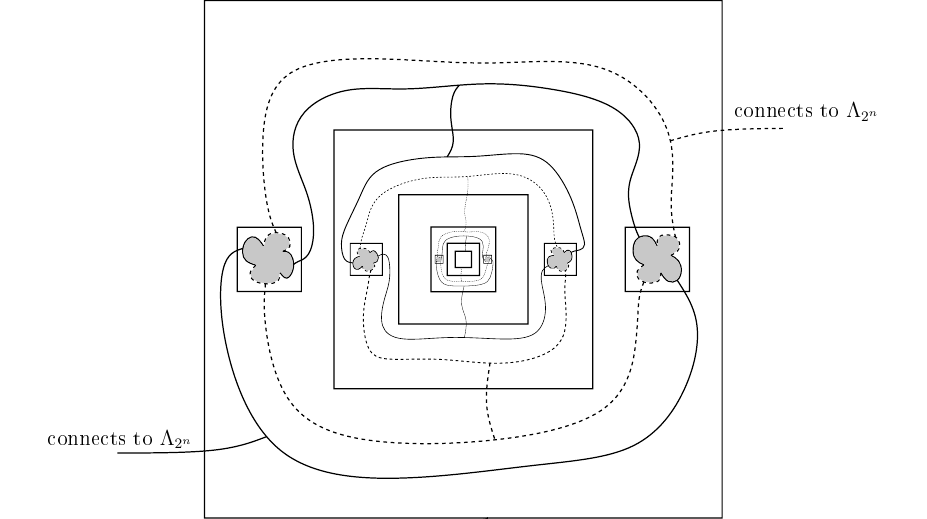}
\caption{A configuration with some good scales in series. The central square is~$\La_{2^m}$ and each annulus represent a different scale. The grey region is the complementary of the explored region~$E$. It is the union of all the flower domains lying within a good scale.} \label{fig:density_good_events}
\end{figure}

Conversely, to have~$\tilde\omega \in A_1(2^m,2^n)$, it suffices for at least one of~$\calF_{k_j}^1$ and~$\calF_{k_j}^2$ to be primally crossed for each~$j$. 
Thus we find 
\begin{align}\label{eq:YoA}
\bbP\big[\tilde\omega\in A_1(2^m,2^n)\,\big|\,\calE\big] 
= \bbP\big[Y_j^1  + Y_j^2   \geq 1 \text{ for all~$j$} \,\big|\,\calE\big] 
= \prod_{j=1}^\ell \bbP\big[Y_j^1  + Y_j^2   \geq 1 \,\big|\,\calE\big].  
\end{align}
The same holds for~$\omega$, so 
\begin{align}\label{eq:XoA}\nonumber
\bbP[\omega\in A_1(2^m,2^n)\,|\,\calE] 
%\bbP[X_j^1  + X_j^2   \geq 1 \text{ for all~$j$} \,|\,\calE] 
&= \prod_{j=1}^\ell \bbP[X_j^1  + X_j^2   \geq 1\,|\,\calE]\\
&= \prod_{j=1}^\ell \bbP[Y_j^1  + Y_j^2\geq 1\,|\,Y_j^1  + Y_j^2  \leq 1 \text{ and }\calE].
\end{align}
Moreover, using that the couple~$(Y_j^1 ,Y_j^2)$ takes each of the four possibilities~$\{0,1\}^2$ with uniformly positive probability, we get 
\begin{align}
	\bbP\big[Y_j^1  + Y_j^2\geq 1\,\big|\,Y_j^1  + Y_j^2  \leq 1 \text{ and }\calE\big] \leq \tfrac{1}{C}\, \bbP\big[Y_j^1  + Y_j^2\geq 1\,\big|\,\calE\big],
\end{align}
for some constant~$C > 1$ independent of~$j$ and the realisations of $\omega$ and $\omega'$ on $E$. 
Inserting this into \eqref{eq:YoA} and \eqref{eq:XoA} we conclude that 
\begin{align}
\bbP\big[\omega\in A_1(2^m,2^n)\,\big|\,\calE\big] 
\leq \big(\tfrac{1}{C}\big)^\ell\, \bbP\big[\tilde\omega\in A_1(2^m,2^n)\,\big|\,\calE\big] 
\leq \big(\tfrac{1}{C}\big)^\ell\, \bbP\big[\omega'\in A_1(2^m,2^n)\,\big|\,\calE\big],
\end{align}
with the second inequality due to the increasing nature of the coupling between~$\omega$,~$\tilde \omega$ and~$\omega'$.
Averaging over all~realisations of $(\omega,\omega')$ on $E$ with~$\omega \in D(K)$ we conclude that 
\begin{align}
\bbP\big[\omega\in A_1(2^m,2^n) \text{ and }  \omega \in D(K)\big]
\leq \big(\tfrac{1}{C}\big)^K \,\bbP\big[\omega'\in A_1(2^m,2^n) \text{ and }  \omega \in D(K)\big].
\end{align}
Divide by~$\bbP[ \omega \in D(K)]$ to obtain the desired conclusion.
\end{proof}

\section{Proof of the main result}\label{section:proof_main}
\begin{proof}[Proof of Theorem~\ref{thm:comparison_exponents}]
Fix~$\delta > 0$ given by Lemma~\ref{lem:density_good_scales} and~$C > 1$ given by Lemma~\ref{lem:comparison_under_good_event}. 
For~$m \leq n$, applying Lemma~\ref{lem:comparison_under_good_event} to~$K = \delta (n-m)$ we find 
\begin{align*}
\phi[A_{1}(2^m,2^n)]&= \bbP[\omega' \in A_{1}(2^m,2^n)]\\
&\geq \bbP[\omega \in D_{m,n}(K) ~\mbox{and}~ \omega' \in A_{1}(2^m,2^n)]\\
&\geq  C^{\delta (n-m)} \, \bbP[\omega \in D_{m,n}(K) ~\mbox{and}~ \omega \in A_{1}(2^m,2^n)]\\
&=  C^{\delta (n-m)}  \,\bbP[\omega \in D_{m,n}(K)\,|\,\omega \in A_{1}(2^m,2^n)]\bbP[\omega \in A_{1}(2^m,2^n)]\\
&=  C^{\delta (n-m)}  \, \phi[D_{m,n}(K) \,|\, A_{01}(2^m,2^n)]\phi[A_{01}(2^m,2^n) \,|\, A_{0}(2^m,2^n)]\\
&\geq \delta C^{\delta (n-m)}  \, \phi[  A_{01}(2^m,2^n) \,|\, A_{0}(2^m,2^n)],
%\delta \bbP[\omega \in A_{01}(2^m,2^n)] &\leq \bbP[\omega \in D_n(K) \cap A_{01}(2^m,2^n)]
\end{align*}
where the second inequality uses Lemma~\ref{lem:comparison_under_good_event} and the last Lemma~\ref{lem:density_good_scales}. Multiply by~$\phi[A_{0}(2^m,2^n)]$ to obtain \eqref{ineq:improved_FKG} for~$r = 2^{m}$ and~$R = 2^{n}$. 
The proof adapts readily to the case where~$r$ and~$R$ are not integer powers of~$2$. 
\end{proof}

%%%%%%%%%%%%%%%%%%%%% Bibliographie %%%%%%%%%%%%%%%%%%%

\bibliographystyle{amsplain}
\bibliography{biblio.bib}

\end{document}